\documentclass[10pt,twoside]{article}
\usepackage{Latex-document}

\newtheorem{thm}{Theorem}

\newtheorem{rem}{Remark}

\newtheorem{question}{Question}
\newcommand {\Bbb}{}

\title{\bf On Some Conformally Invariant Fully\vskip -2mm Nonlinear Equations}
\author{YanYan Li\thanks{Department of Mathematics, Rutgers University, 110 Frelinghuysen Rd.,
Piscataway, NJ 08854, USA. E-mail:
yyli@math.rutgers.edu}\vspace*{-0.5cm}}
\date{\vspace{-8mm}}

\begin{document}

\maketitle

\thispagestyle{first} \setcounter{page}{177}

\begin{abstract}

\vskip 3mm

We will report some results concerning  the Yamabe problem  and the Nirenberg problem.  Related topics will also
be discussed. Such studies  have  led to  new results on some conformally invariant fully nonlinear equations
arising from geometry. We will also present these results which include some Liouville type theorems, Harnack type
inequalities, existence and compactness of solutions to some nonlinear version of the Yamabe problem.

\vskip 4.5mm


\end{abstract}

\vskip 12mm

In this talk, we present some recent joint work with Aobing Li \cite{LL} on some conformally invariant fully
nonlinear equations.

For $n\ge 3$, consider
\begin{equation}
-\Delta u=\frac {n-2}2 u^{ \frac{n+2}{n-2}},
\qquad \mbox{on} \quad \Bbb {R}^n.
\label{eq1new}
\end{equation}

The celebrated Liouville type theorem of Caffarelli, Gidas and Spruck
(\cite{CGS}) asserts that positive $C^2$ solutions of
(\ref{eq1new}) are of the form
$$
u(x)=(2n)^{ \frac {n-2}4 }
\left(\frac a {1+a^2|x-\bar{x}|^2}\right)^{\frac{n-2}{2}},
$$
where $a>0$ and $\bar x\in \Bbb R^n$.
Under an additional decay hypothesis
$u(x)=O(|x|^{2-n})$, the result was proved
by Obata (\cite{O}) and Gidas, Ni and Nirenberg (\cite{GNN}).

Let $\psi$ be a M\"obius transformation.
$$
\left(  u_\psi
^{-  \frac {n+2}{n-2} }\Delta u_\psi
\right)=
\left(  u^{  -\frac {n+2}{n-2} }\Delta u
\right)\circ\psi,
\qquad\mbox{on}\ \Bbb R^n,
$$
where  $u_\psi:=|J_\psi|^{ \frac {n-2}{2n} }(u\circ \psi)$
and  $J_\psi$ denotes the Jacobian of $\psi$.
In particular, if  $u$ is a positive solution of
(\ref{eq1new}), so is $u_\psi$.

We call a fully nonlinear operator $H(x,u, \nabla u, \nabla^2 u)$
conformally invariant on $\Bbb R^n$ if
for any M\"obius transformation
$\psi$ and any positive function $u\in C^2(\Bbb R^n)$
\begin{equation}
H(\cdot, u_\psi, \nabla u_\psi, \nabla^2 u_\psi)\equiv
 H(\cdot, u, \nabla u, \nabla^2 u)\circ \psi.
\label{conformal3}
\end{equation}
We showed in \cite{LL} that $H(x,u, \nabla u, \nabla^2 u)$
is conformally invariant if and only if
$$
H(x,u, \nabla u, \nabla^2 u)\equiv F(A^u),
$$
where
\begin{equation}
A^u:= -\frac{2}{n-2}u^{  -\frac {n+2}{n-2} }
\nabla^2u+ \frac{2n}{(n-2)^2}u^ { -\frac {2n}{n-2} }
\nabla u\otimes\nabla u-\frac{2}{(n-2)^2} u^ { -\frac {2n}{n-2} }
|\nabla u|^2I,
\label{2}
\end{equation}
and  $F$ is invariant under orthogonal conjugations.

Let $U$ be an open subset of $n\times n$
symmetric matrices which is  invariant under orthogonal conjugations
(i.e.
$
O^{-1}UO=U
$ for all orthogonal matrices $O$)
and has the property that $U\cap \{M+tN\ |\ 0<t<\infty\}$ is convex
for any $n\times n$ symmetric matrix $M$ and any
$n\times n$ positive definite symmetric matrix $N$.

Let $F\in C^1(U)$ be  invariant under orthogonal conjugation and
be elliptic, i.e.
$$
\left(F_{ij}(M)\right)>0,\qquad \forall\ M\in U,
$$
where $F_{ij}(M):=\frac{\partial F}{ \partial M_{ij} }(M)$.

The following theorem extends the result of Obata and Gidas, Ni and Nirenberg
to all conformally invariant operators of elliptic type.
\begin{thm} {\rm (\cite{LL})} For $n\ge 3$, let $U$ and $F$ be
as above, and
let $u\in C^2(\Bbb R^n)$ be a positive solution of
$$
F(A^u)=1, \qquad \mbox{on}\   \Bbb R^n.
$$
Assume that $u$ is regular at infinity, i.e.,
$|x|^{2-n}u(x/|x|^2)$ can be
extended to a
positive $C^2$ function  near the origin.
Then
for some $\bar x\in \Bbb R^n$ and for  some positive constants $a$ and $b$,
$$
u(x)\equiv \left(\frac a {1+b^2|x-\bar{x}|^2}\right)^{\frac{n-2}{2}},
\qquad \forall\ x\in \Bbb R^n.
$$
\label{thm4-1}
\end{thm}

\begin{rem}
{\rm  In fact, as established in \cite{LL}, the conclusion of the
above theorem still holds when replacing the assumption $u\in
C^2(\Bbb R^n)$ by a weaker assumption that $u\in C^2(\Bbb
R^n\setminus\{0\})$,  $u$ can be extended
 to a
positive continuous  function  near the origin, and $ \lim_{x\to
0}(|x||\nabla u(x)|)=0$.} \label{remark1}
\end{rem}

Theorem \ref{thm4-1} indicates that behavior of solutions
to conformally invariant equations is very rigid.
Thus we expect some good theories for
conformally invariant uniformly elliptic fully nonlinear
equations.   Let $F$ be $C^\infty$ functions defined
on $n\times n$ real symmetric matrices, and let
$F$ be invariant under orthogonal conjugations.  We  assume
that for some constants $0<\lambda\le \Lambda<\infty$,
$$
\lambda I\le \left(F_{ij}(M)\right)\le \Lambda I,
\qquad \mbox{for all}\ n\times n\
\mbox{real symmetric matrices}.
$$
We raise the following

\begin{question}  {\rm Let $F$ be as above,
and let $B_1$ be a unit ball in $\Bbb R^n$ and $a>0$ be some constant.
Are there some positive constants $\alpha$ and $C$, depending only
on $F$, $a$ and $n$ such that
for any  positive $C^\infty$ solution
 $u$
 of
$$
F\left(A^u\right)=0, \qquad \mbox{in}\  B_1
$$
satisfying
$$
\min_{ \overline B_1}u\ge a, \qquad \|u\|_{ C^2(B_1) }\le \frac 1a,
$$
we have
$$
\|u\|_{ C^{2,\alpha}(B_{\frac 12})  }\le C?
$$}
\label{question1}
\end{question}

Other interesting questions include to understand
behavior near an isolated singularities of a solution
in a punctured disc of  this subclass of uniformly elliptic
equations and to  establish
some removable singularity results.

Let $(M,g)$ be an $n-$dimensional smooth
Riemannian manifold without boundary, consider the
 Schouten
 tensor
$$
A_g=\frac{1}{n-2}\left(Ric_g-\frac{R_g}{2(n-1)}g\right),
$$
where $Ric_g$ and $R_g$ denote respectively the
Ricci tensor and the scalar curvature associated with $g$.

For  $1\le k\le n$,
let
$$
\sigma_k(\lambda)=\sum_{1\le i_1<\cdots<i_k\le n}\lambda_{i_1}
\cdots \lambda_{i_k},\qquad\qquad \lambda=(\lambda_1,
 \cdots, \lambda_n)\in \Bbb R^n,
$$
denote the $k-$th symmetric function, and let
 $\Gamma_k$ denote the
connected component of $\{\lambda\in \Bbb R^n\ |\
\sigma_k(\lambda)>0\}$ containing  the
positive cone $\{\lambda\in \Bbb R^n\ |\
\lambda_1, \cdots, \lambda_n>0\}$.

It is known (see, e.g., \cite{CNS}) that $\Gamma_k$ is
a convex cone with its vertex at the origin,
$$
\Gamma_n\subset \cdots\subset \Gamma_2\subset \Gamma_1,
$$
$$
\frac{\partial \sigma_k}{\partial \lambda_i}>0\quad
\mbox{in}\ \Gamma_k, \ 1\le i\le n,
$$
and
$$
\sigma_k^{\frac 1k}\ \mbox{is concave in}\ \Gamma_k.
$$

Fully nonlinear elliptic equations involving
$\sigma_k(D^2u)$ have been investigated in the
classical and pioneering paper of
 Caffarelli, Gidas and Nirenberg
\cite{CNS}.  For extensive studies and outstanding
results on such equations, see, e.g.,
Guan and Spruck \cite{GS}, Trudinger \cite{T1},
Trudinger and Wang \cite{TW}, and the references therein.
On Riemannian manifolds of nonnegative curvature, Li  studied
in \cite{Li} equations
\begin{equation}
\sigma_k^{\frac 1k}\left(\lambda(\nabla^2_gu+g)\right)=\psi(x,u),
\label{neg}
\end{equation}
where $\lambda(\nabla^2_gu+g)$ denotes eigenvalues of
$\nabla^2_gu+g$ with respect to
$g$.
On general Riemannian manifolds, Viaclovsky
introduced and systematically studied
in \cite{V3} and \cite{V1} equations
\begin{equation}
\sigma_k^{\frac 1k}\left(\lambda(A_g)\right)=\psi(x,u),
\label{c1new}
\end{equation}
where $\lambda(A_g)$
 denotes the eigenvalues of $A_g$ with respect to $g$.
On $4-$dimensional general Riemannian manifolds,
remarkable results on (\ref{c1new}) for $k=2$
were
 obtained by Chang, Gursky and Yang in
\cite{CGY1} and \cite{CGY2}, which include Liouville type theorems,
existence and compactness of solutions,
as well as applications to topology.
On the other hand, works on the Yamabe equation by
Caffarelli, Gidas and Spruck (\cite{CGS}),
Schoen (\cite{S1} and \cite{S2}), Li and Zhu (\cite{LZhu}), and
Li and Zhang (\cite{LZ}),  have played an important
role in our approach to the study of
(\ref{c1new}) as developed in \cite{LL}.

Consider
\begin{equation}
\sigma_k\left(\lambda(A_g)\right)=1,
\label{c1}
\end{equation}
 together with
\begin{equation}
\lambda(A_g)\in \Gamma_k.
\label{c2}
\end{equation}
Let $g_1=u^{\frac 4 {n-2} }g_0$ be a conformal change of
metrics, then (see, e.g., \cite{V3}),
$$
A_{g_1}=-\frac{2}{n-2}  u^{-1}\nabla_{g_0}^2u+ \frac{2n}{(n-2)^2} u^{-2}
\nabla_{g_0}u \otimes\nabla_{g_0}u -\frac{2}{(n-2)^2} u^{-2}
|\nabla_{g_0}u|_{g_0}^2g_0+A_{g_0}.
$$

Let $g=u^{\frac{4}{n-2}}g_{flat}$, where $g_{flat}$ denotes
the Euclidean metric on $\Bbb R^n$.  Then by the above transformation formula,
$$
A_g=u^{\frac{4}{n-2}}A^u_{ij}dx^idx^j,
$$
where $A^u$ is given by (\ref{2}).

Equations (\ref{c1}) and (\ref{c2}) take the form
\begin{equation}
\sigma_k(\lambda(A^u))=1,\qquad \mbox{on} \quad \Bbb {R}^n,
\label{eq1}
\end{equation}
and
\begin{equation}
\lambda(A^u)\in \Gamma_k,\qquad \mbox{on} \quad  \Bbb R^n.
\label{gammacone}
\end{equation}

Our next result extends the Liouville
type theorem of Caffarelli, Gidas and Spruck to
 all $\sigma_k$, $1\le k\le n$.  For $k=1$,
equation
(\ref{eq1}) is
(\ref{eq1new}).

\begin{thm} {\rm (\cite{LL})}\
  For $n\ge 3$ and $1\le k\le n$, let $u\in C^2(\Bbb R^n)$
be
a positive solution of (\ref{eq1}) satisfying (\ref{gammacone}).
Then for some $a>0$ and $\bar x\in \Bbb R^n$,
\begin{equation}
u(x)=c(n,k)\left(\frac a {1+a^2|x-\bar{x}|^2}\right)^{\frac{n-2}{2}},
\qquad \forall\ x\in \Bbb R^n,
\label{t2e}
\end{equation}
where $c(n,k)=2^{ (n-2)/4 }{n \choose
 k}^{(n-2)/4k}$.
\label{t2}
\end{thm}

  The case $k=2$ and $n=4$
was obtained by Chang, Gursky and Yang (\cite{CGY2}).
More recently, they (\cite{CGY3})
have independently established the result for $k=2$ and $n=5$, and they
also established the result for $k=2$ and $n\ge 6$ under
an additional hypothesis $\int_{\Bbb R^n}u^{ \frac {2n}{n-2} }<\infty$.
Under an  additional hypothesis that
$\frac 1{|x|^{n-2}}u(\frac x{|x|^2})$ can be extended
to a $C^2$ positive  function near $x=0$,
the case $2\le k\le n$ was obtained by
Viaclovsky (\cite{V3}, \cite{V2}).
 As mentioned above, the case $k=1$ was obtained
by Caffarelli, Gidas and Spruck, while under
an additional hypothesis that
$\frac 1{|x|^{n-2}}u(\frac x{|x|^2})$ is bounded
near  $x=0$,
the case $k=1$ was obtained by Obata, and by Gidas, Ni and Nirenberg,

The methods of Chang, Gursky and Yang in \cite{CGY2} and \cite{CGY3}
include an ingenious way of using the Obata technique which, as
 they pointed out, allows the possibility to be generalized to
establish the uniqueness of solutions on general
Einstein manifolds.
Our proof of  Theorem \ref{t2} is very different
from that of \cite{CGY2} and \cite{CGY3}.
 A
crucial ingredient in our proof
is the following Harnack type inequality.

\begin{thm} {\rm (\cite{LL})}\
For $n\ge3$, $1\le k\le n$, and $R>0$, let
$B_{3R}\subset \Bbb R^n$ be a ball of radius $3R$ and let
 $u\in C^2(B_{3R})$ be a positive solution of
\begin{equation}
\sigma_k(A^u)=1, \qquad \mbox{in}\ B_{3R},
\label{yamabe}
\end{equation}
satisfying
\begin{equation}
\lambda(A^u)\in \Gamma_k,\qquad \mbox{in}\  B_{3R}.
\label{condition}
\end{equation}
Then
\begin{equation}
(\max_{\overline{B}_R} u)(\min_{\overline{B}_{2R}}u) \le C(n)R^{2-n}.
\label{4}
\end{equation}
\label{t1}
\end{thm}

 The above Harnack type inequality for
$k=1$ was obtained by Schoen (\cite{S2}) based
on the Liouville type theorem of Caffarelli, Gidas and Spruck.  An important
step toward our proof of Theorem \ref{t1} was
taken in an earlier work of Li and
Zhang (\cite{LZ}), where they gave a different proof
of Schoen$'$s   Harnack type inequality without using the
Liouville type theorem.

Our next result concerns existence and compactness of solutions.
\begin{thm} {\rm (\cite{LL})}\ For $n\ge 3$ and $1\le k\le n$, let $(M, g)$ be an
$n-$dimensional  smooth  compact locally conformally flat
 Riemannian  manifold without boundary satisfying
$$
\lambda\left(A_{g}\right)\in \Gamma_k,\qquad \mbox{on}\  M.
$$
 Then there exists some smooth positive function $u$ on $M$ such that
$\hat g=u^{\frac 4{n-2}} g$ satisfies
\begin{equation}
\lambda\left(A_{\hat g}\right)\in \Gamma_k,\ \
\sigma_k\left(\lambda(A_{\hat g})\right)=1,\qquad \mbox{on}\ M.
\label{ghat}
\end{equation}
Moreover, if  $(M,g)$ is not conformally diffeomorphic
to the  standard $n-$sphere,  all solutions of the above
satisfy, for all $m\ge 0$, that
$$
\|u\|_{C^m(M,g)}+\|u^{-1}\|_{C^m(M,g)}\le C,
$$
where $C$ depends only on $(M,g)$ and $m$.
\label{ll1}
\end{thm}

For $k=1$, it is the Yamabe
problem for  locally conformally flat manifolds
with positive Yamabe invariants,
 and  the result  is due to Schoen (\cite{S0}-\cite{S1}).
The Yamabe problem was solved through the work of
Yamabe \cite{Y}, Trudinger \cite{Tr}, Aubin \cite{A}, and Schoen \cite{S0}.
For $k=2$ and $n=4$,  the result was proved without the
locally conformally flatness hypothesis by
Chang, Gursky and Yang \cite{CGY2}. For $k=n$, the existence
result was established by Viaclovsky \cite{V1} for a class of manifolds
which are not necessarily locally conformally flat.
For $k\neq \frac n2$, the result is independently obtained by
Guan and Wang in \cite{GW2} using a heat flow method.
More recently, Guan,  Viaclovsky and
Wang \cite{GVW} have proved
that  $\lambda(A_g)\in \Gamma_k$ for
$k\ge \frac n2$ implies  the positivity of the  Ricci tensor,
and therefore, by classical results,
$(M,g)$ is conformally covered by $\Bbb S^n$ and the existence
and compactness  results in this case follow easily.

Our proof of  Theorem  \ref{thm4-1}, different from
the ones  in \cite{O}, \cite{GNN},
\cite{CGS}, \cite{V3} and \cite{V2}, is in the spirit
of the new proof of the Liouville type theorem of
Caffarelli, Gidas and Spruck given by Li and Zhu in
\cite{LZhu}.
We  also make use of the substantial
simplifications of Li and Zhang in \cite{LZ}
to  the proof in \cite{LZhu}.
 The proof is along the line of the
pioneering work of  Gidas, Ni and Nirenberg
\cite{GNN}, which in particular does not need the kind of divergence
structure needed for
the method of Obata \cite{O} and therefore can be applied
in much more generality.

In our proofs blow up arguments are used, which require
local derivative estimates of solutions.  For
$\sigma_1$ (the Yamabe equation), such  estimates
follow from
standard elliptic theories.
Guan and Wang \cite{GW} established
 local gradient and second
derivative estimates for  $\sigma_k$, $k\ge 2$.
 Global
gradient and second derivative estimates
for $\sigma_k$  were
obtained by Viaclovsky
\cite{V1}.  For the related equation (\ref{neg}) on manifolds of nonnegative
curvature,
global gradient and second derivative estimates
 were obtained
by Li  in \cite{Li}.
By  the concavity of $\sigma_k^{\frac 1k}$,
$C^{2,\alpha}$ estimates hold due to the classical work of
 Evans \cite{E}
and Krylov \cite{K}.
For the proof of the existence part of
Theorem \ref{ll1}, we introduce a homotopy
$\sigma_k(t\lambda+(1-t)\sigma_1(\lambda))$,
defined on $(\Gamma_k)_t=\{
\lambda\in \Bbb R^n\ |\ t\lambda+(1-t)\sigma_1(\lambda)\in
\Gamma_k\}$, which establishes a natural link
between (\ref{ghat}) and  the Yamabe  problem.
We extend the  local estimates in \cite{GW}
 for $\sigma_k$ to
$\sigma_k(t\lambda+(1-t)\sigma_1(\lambda))$,
with estimates uniform in $0\le t\le 1$.
The compactness results as stated
in Theorem \ref{ll1} were  established in \cite{LL}
along the homotopy.  The  compactness results
for the Yamabe problem was established by Schoen \cite{S1}.
We gave a different proof which does
not rely on the Liouville type theorem, which allows us
to establish existence results for more general
$f$ than $\sigma_k$ for which Liouville type theorems
are not available.
The existence results follow from the compactness results
with the help of the degree theory
for second order fully nonlinear elliptic operators
(\cite{Li1})
as well as the degree counting formula
for the Yamabe  problem (\cite{S1}).

The first step in our proof of the Liouville type Theorem \ref{t2}
is to establish the Harnack type
inequality (Theorem \ref{t1}), from which we obtain
 sharp asymptotic behavior at infinity
of an entire solution.    Then we establish  Theorem~\ref{t2}
by distinguishing into two cases.  In the case
$k>\frac n2$,  Theorem~\ref{t2} is proved by using
the  sharp asymptotic behavior of
an entire solution and Theorem \ref{thm4-1}-Remark \ref{remark1},
 together with a result of Trudinger and Wang
(\cite{TW}).
In the case $1\le k\le \frac n2$,  Theorem~\ref{t2} is proved by
the  sharp asymptotic behavior of
an entire solution together with the Obata type integral formula
of Viaclovsky (\cite{V3}).  For the second case, divergence
structure of the equation is used.

Theorem \ref{t2}, Theorem \ref{t1} and Theorem \ref{ll1}
are established for more general nonlinear $f$ than $\sigma_k$
in \cite{LL}, including those
for which  no divergence structure is available.

\label{lastpage}

\end{document}